\documentclass[preprint,10pt,nopreprintline]{elsarticle}

\usepackage[english]{babel}
\usepackage{geometry}
\usepackage{amsmath,amsfonts,amssymb,amsthm}
\usepackage[final]{hyperref}
\usepackage{graphicx}
\usepackage{empheq}
\usepackage{mathtools}
\usepackage{physics}
\usepackage{mathtools}
\usepackage[dvipsnames]{xcolor}
\usepackage{pdflscape}
\usepackage{etoolbox}
\usepackage{fancyhdr}
\usepackage{lastpage}
\usepackage{ifdraft}
\usepackage{hyphenat}
\usepackage{natbib}
\ifdraft{\usepackage{showlabels}}{}
\usepackage{enumitem}
\usepackage{siunitx}
\usepackage{booktabs}
\usepackage{multirow}
\usepackage{array}
\usepackage{caption}
\usepackage{float}
\usepackage{rotating}
\usepackage[noabbrev]{cleveref}
\usepackage{algorithm}
\usepackage[noend]{algpseudocode}

\usepackage{lineno}

\begin{document}

\begin{frontmatter}



\title{Exploring Memory Effects: Sparse Identification in Vector-Borne Diseases} 


\author[1]{Dimitri Breda}
\ead{dimitri.breda@uniud.it}
\author[1]{Muhammad Tanveer}
\ead{tanveer.muhammad@spes.uniud.it}
\author[2]{Jianhong Wu}
\ead{wujh@yorku.ca}
\author[3]{Xue Zhang}
\ead{zhangxue@mail.neu.edu.cn}

\affiliation[1]{organization={CDLab – Computational Dynamics Laboratory Department of Mathematics, Computer Science and Physics – University of Udine},
                addressline={via delle Scienze 206}, 
                city={Udine},
                postcode={33100}, 
                country={Italy}}
\affiliation[2]{organization={Laboratory for Industrial and Applied Mathematics Department of Mathematics and Statistics -- York University},
                addressline={Toronto, ON}, 
                city={Toronto},
                postcode={M3J 1P3}, 
                country={Canada}}
\affiliation[3]{organization={Department of Mathematics, Northeastern University},
                addressline={Shenyang}, 
                city={Shenyang},
                postcode={110819}, 
                country={China}}

\begin{abstract}
Predicting the human burden of vector-borne diseases from limited surveillance data remains a major challenge, particularly in the presence of nonlinear transmission dynamics and delayed effects arising from vector ecology and human behavior. We develop a data-driven framework based on an extension of Sparse Identification of Nonlinear Dynamics (SINDy) to systems with distributed memory, enabling discovery of transmission mechanisms directly from time series data. Using severe fever with thrombocytopenia syndrome (SFTS) as a case study, we show that this approach can uncover key features of tick-borne disease dynamics using only human incidence and local temperature data, without imposing predefined assumptions on human case reporting. We further demonstrate that predictive performance is substantially enhanced when the data-driven model is coupled with mechanistic representations of tick–host transmission pathways informed by empirical studies. The framework supports systematic sensitivity analysis of memory kernels and behavioral parameters, identifying those most influential for prediction accuracy. Although the approach prioritizes predictive accuracy over mechanistic transparency, it yields sparse, interpretable integral representations suitable for epidemiological forecasting. This hybrid methodology provides a scalable strategy for forecasting vector-borne disease risk and informing public health decision-making under data limitations.   
\end{abstract}




\end{frontmatter}



\section{Introduction}
Projecting tick-borne infection risk to human under substantial uncertainties of vector abundance, environmental condition and human-tick interaction imposes significant computational challenges since the data we observe, i.e., the human incidence of the relevant diseases, is only the final outcome of a complex and multi-year tick-host pathogen transmission and tick-human contact processes. Evidence-based short-term forecasting requires rigorous computational and data fitting techniques to estimate parameters of the underlying transmission dynamics model and long-term scenario analysis (for example, under projected climate scenarios) demands identifying the first-principle mechanisms that govern the patterns of tick-host and host-human contact for the human case reported. The challenge is amplified by the limited availability of time series of human incidence and environmental conditions despite a large body of fields and laboratory studies on ticks and hosts.

To address these short term forecasting and long term projection challenges, we employ the Sparse Identification of Nonlinear Dynamics (SINDy)~\cite{bpk16} framework that automatically discovers governing equations in mechanistic models from time series data of some state variables. SINDy was pioneered for mechanistic models using ordinary differential equations (ODEs)~\cite{bk19,brunton2022data,bpk16}, and was later extended to those models involving spatial variation such as partial differential equations (PDEs)~\cite{rud17}, stochastic variation such as stochastic differential equations~\cite{boninsegna2018sparse}, and more recently, development and feedback delay such as differential equations with constant (discrete) time delays~\cite{bbt24,kopeczi23,pec24,sandoz2023sindy}. These advancements highlight the flexibility and applicability of the framework, reinforcing its significance as a pivotal technique in the rapidly expanding field of data-driven scientific discoveries.

Despite this substantial progress, a significant gap that is highly relevant to projection of vector-borne disease burden in human persists in the application of SINDy to systems characterized by high nonlinearity and distributed memory effects. Distributed delays (DDEs) and renewal-type Volterra integral equations (REs) involving integration over a continuum of past states to model distributed memory dynamics are key ingredients in modelling biological, ecological and engineering systems in which the evolution of the system relies on a range of past states~\citep{beretta2016discrete,emo_renewal21,metz2014dynamics}. Examples include age-structured population models~\cite{metz2014dynamics}, maturation processes in cell populations~\cite{emo_renewal21}, tick-borne disease transmission dynamics \cite{CSIAM-LS-1-2} and genetic regulatory networks with complex time lags \cite{glass2021nonlinear}.

Building on the recently developed SINDy extension for distributed memory \cite{bredadistributed2025}, we show how effective the SINDy framework is for nearcasting and short-term forecasting tick-borne disease transmission risk in a temporally-changing environment. Tick-borne diseases present an ideal test case: they involve distributed memory effects across multiple biological timescales from tick maturation and pathogen incubation periods to human immune response dynamics. Traditional compartmental models \cite{wu2020transmission,CSIAM-LS-1-2} struggle to capture these continuous, distributed processes without explicitly specifying many intermediate biological stages and transition rates. By contrast, our data-driven approach automatically discovers these distributed delay structures directly from human incidence data, without requiring prior specification of tick life stages, vector-pathogen dynamics, or disease reporting mechanisms. Our method combines integral-based formulations~\citep{messenger2021weak,schaeffer2017sparse,wei2022sparse} that improve noise robustness with structured candidate library construction that explicitly includes distributed delay variables and quadrature approximations for memory kernel integration.

While our framework is readily transferable to other vector-borne systems characterized by nonlinear interactions and delayed transmission, we take Severe Fever with Thrombocytopenia Syndrome (SFTS) as an example to demonstrate SINDy methodology, applying it to the most affected city Dalian of the Liaoning Province, China. SFTS is an emerging tick-borne infectious disease caused by the Bunyavirus/SFTS virus (SFTSV) and transmitted primarily by {\it Haemaphysalis longicornis} \cite{kim2013severe}. With fatality rates ranging from 12–50\% in East Asia \cite{he2020severe}, SFTS has been recognized by the World Health Organization as a priority emerging zoonosis since its first identification in China in 2009. The disease remains severe in central China, has been expanding to higher latitudes and has recently emerged in northeastern provinces such as Liaoning, where incidence shows a clear upward trend. Historically, nearly 90\% of reported SFTS cases in Liaoning since 2011 were concentrated in Dalian and Dandong. However, recent epidemiological trends reveal a significant shift: during the period from 2022 to 2024, Dalian experienced a more rapid increase in incidence, accounting for over 60\% of provincial cases. This establishes Dalian as the current epicenter of SFTS within the province of Liaoning. Given this background of intensifying and localized transmission, we apply our data-driven method to SFTS case data from Dalian. Our method directly optimizes for predictive accuracy and identifies distributed delay effects without requiring explicit specification of intermediate biological processes.


The remainder of this work is organized as follows: Section~\ref{sec:dde_re} reviews the mathematical foundation of REs and our adopted method for identifying distributed delay systems. This section also details the construction of a library based on quadrature and the formulation of sparse regression. Section~\ref{sec:sfts} provides further background discussions of SFTS and describes the epidemiological data used in the case study.  Section~\ref{sec:sfts_sindy} presents numerical experiments demonstrating the performance of our method on SFTS transmission data, focusing on examining both model identification accuracy and parameter sensitivity. Finally, we conclude Section~\ref{sec:conclusion} with a discussion of the wider implications of this study and outline future research directions for data-driven discoveries in infinite-dimensional systems arising from biological systems.


\section{SINDy for distributed delay renewal equations} \label{sec:dde_re}
Many dynamical systems in practice depend not only on discrete delays but also on the cumulative history of the state. Rather than a system responding to its value at one or more specific time points, the current behavior and behavioral change reflects a weighted average of how the system evolved. This distributed memory where the present evolves from contributions across the entire past history is fundamentally different from systems with discrete, isolated delays and leads naturally to REs and related integral formulations. REs and their generalizations through PDEs (for those interested in models of structured populations, we suggest \cite{Boldin02082024,franco2023modelling} for further reading) provide powerful frameworks for modeling such systems, particularly in structured population dynamics.

The core structure of a RE captures this distributed dependence through a functional relationship as
\begin{equation}\label{eq:re}
    x(t) = G(x_t)
\end{equation}
which connects the current state directly to its history $x_t$ defined as $x_t(s) \coloneqq x(t+s)$ for $s \in [-\tau, 0]$, where $\tau > 0$ is the maximum delay. Here, $G: X \to \mathbb{R}^n$ is a functional operating on a Banach space $X$ of history functions mapping $[-\tau,0]$ to $\mathbb{R}^n$ with appropriate regularity (typically measurable for REs \cite{dgg07}). The functional $G$ accounts for how the system aggregates influence from its past.
In many applications,  $G$ is defined by an integral
\begin{equation}\label{eq:G}
    G(\varphi)\coloneqq \int_{-\tau}^{0} g(s, \varphi(s))\dd s,\quad\varphi\in X,
\end{equation}
where the kernel function $g: [-\tau,0] \times \mathbb{R}^n \to \mathbb{R}^n$ is a typically nonlinear map that assigns weights to past states. The kernel $g$ is the decisive element; it determines the character of the model by specifying which parts of the history matter (first argument of $g$)  and how much they matter at each state value (via the second argument of $g$).
As we clarify, this nonautonomous dependence of $g$ on $s$ is a key aspect for its reconstruction from observed data.


Identifying the kernel $g$, the first principle, from observed data is the central challenge in reconstructing REs from measurements. Following SINDy framework~\citep{messenger2021weak,schaeffer2017sparse, wei2022sparse}, here we recall from  \cite{bredadistributed2025} a recently introduced framework based on quadrature approximations that allows us to adapt SINDy methodology to discover the right-hand side \eqref{eq:G} of an RE \eqref{eq:re} .
In practical applications, states are observed at discrete times, often accompanied by measurement noise. The available data consist of $m$ snapshots organized in a data matrix as follows
\begin{equation}\label{eq:data}
    \mathbf{X} =({x}_{1},\ldots,{x}_{n})\coloneqq
\begin{pmatrix}
{x}_{1}(t_{1}) & \cdots & {x}_n(t_{1}) \\
\vdots & \ddots & \vdots \\
{x}_{1}(t_m) & \cdots & {x}_n(t_m)
\end{pmatrix}
\in \mathbb{R}^{m \times n}.
\end{equation}


The first step in using SINDy to find the kernel $g$ in REs is to use a suitable quadrature formula to approximate the integral $G(x_{t})$ in \eqref{eq:G}. Thus the integral is transformed into a weighted sum of kernel values:
\begin{equation}\label{eq:quadrature}
\int_{-\tau}^{0} g(s,x(t+s))\dd s\approx\sum_{k=1}^{K} w_{k} g(s_{k},x(t+s_{k})). 
\end{equation}
Here, $K$ denotes the number of quadrature nodes $s_{k} \in [-\tau, 0]$ with associated weights $w_{k}$ for $k = 1, \ldots, K$. The choice of quadrature rule (e.g., rectangle, trapezoidal, or Clenshaw-Curtis \citep{clenshaw1960method,gentleman1972implementing,trefethen2009approximation}) should be based on the smoothness properties of the kernel and the available data resolution.

SINDy framework assumes that the kernel $g$ can be expressed as a sparse linear combination of basis functions from a candidate library $\Theta$. Because $g$ depends on both the delay variable $s$ and the state value $x(t+s)$, the library must be constructed to capture both of these arguments. Importantly, the explicit dependence of the kernel on the delay variable $s$ (nonautonomous dependence) is a key feature that must be preserved during identification, as it determines how memory effects vary across the delay interval.

For each quadrature node $k = 1, \ldots, K$, we construct a library $\Theta_{k}$ of candidate functions that depend on the node $s_{k}$ and the time-shifted state values $x(t+s_{k})$
as

\begin{equation*}
\left\{\setlength\arraycolsep{0.1em}\begin{array}{rcl}
\Theta_{k}\coloneqq\Theta(\sigma_{k},\mathbf{X}(\cdot+\sigma_{k}))&=\big[&1,\\[2mm]
&&\sigma_{k},\mathbf{X}_{1}(\cdot+\sigma_{k}),\ldots,\mathbf{X}_{n}(\cdot+\sigma_{k}),\\[2mm]
&&\sigma_{k}^{2},\sigma_{k}\mathbf{X}_{1}(\cdot+\sigma_{k}),\ldots ,\sigma_{k}\mathbf{X}_{n}(\cdot+\sigma_{k})),\mathbf{X}_{1}^{2}(\cdot+\sigma_{k}),\\[2mm]
&&\ldots\\[2mm]
&&\sigma_{k}^{d},\sigma_{k}^{d-1}\mathbf{X}_{1}(\cdot+\sigma_{k}),\ldots,\mathbf{X}_{n}^{d}(\cdot+\sigma_{k})\big],
\end{array}\right.
\end{equation*}
where polynomial terms are included up to a specified degree $d$. Users can augment this library with additional nonlinear functions based on prior knowledge (e.g., either $s_{k}$, $\mathbf{X}_{j}(\cdot+s_{k})$, $j=1,\ldots,n$, $e^{s_{k}}$, $\sin(\mathbf{X}_{j}(\cdot+s_{k}))$, etc.). The notation $\mathbf{X}_{j}(\cdot+s_{k})$ represents the $j$-th state variable evaluated at time-shifted instants $t_{i}+s_{k}$ for $i=1,\ldots,m$, such shifted values can be recovered through piecewise linear interpolation when necessary.

The complete SINDy library for identifying the REs \eqref{eq:re} is formed by taking a weighted sum of these libraries $\Theta_{k}$, $k=1,\ldots,K$, with weights given by the quadrature weights $w_{k}$. The corresponding identification problem reads
\begin{equation}\label{eq:reg-re}
\mathbf{X}_{j} \approx\left( \sum_{k=1}^{K} w_{k} \Theta(\sigma_{k}, \mathbf{X}(\cdot+\sigma_{k}))\right)\xi_{j}
\end{equation}
for each state component $j = 1, \ldots, n$. In this context, in \eqref{eq:reg-re} the vector $\xi_{j} \in \mathbb{R}^{p}$ represents a sparse vector of coefficients that determines which library terms are actively involved in the kernel function for the $j$-th equation component. Each of these sparse vectors is derived by solving the regularized least-squares regression problem

\begin{equation*}
\begin{split}
\hat{\xi}_{j} = &\arg\min_{\xi\in\mathbb{R}^{p}} \left\| \mathbf{X}_{j} - \left( \sum_{k=1}^{K} w_{k} \Theta(\sigma_{k}, \mathbf{X}(\cdot+\sigma_{k}))\right){\xi}_{j} \right\|_{2}\\
&+\lambda\|{\xi}_{j}\|_{1}.
\end{split}
\end{equation*}
for each $j=1,\ldots,n$. This represents a standard convex optimization problem that can be effectively addressed using tools like STLS \citep{bk19} and LASSO \cite{lasso96}. By assembling all the resulting sparse vectors into a matrix $\Xi\coloneqq (\hat{\xi}_{1},\ldots,\hat{\xi}_n)\in\mathbb{R}^{p\times n}$, we finally identify \eqref{eq:re} with $x(t)\approx\Theta(x^{T}(t))\Xi$, thereby uncovering the fundamental analytical structure of the kernel $g$ in \eqref{eq:G}. This approach allows practitioners to derive interpretable, sparse functional forms that control distributed memory effects in REs and makes explicit the nonautonomous structure of the kernel.
\section{Severe Fever with Thrombocytopenia Syndrome in Dalian}\label{sec:sfts}
SFTS is an emerging disease transmitted from ticks to humans, first identified in central China in 2009. Since then, it has spread to South Korea, Japan, and other East Asian countries. The disease is caused by the SFTS virus (SFTSV), which belongs to the \textit{genus Bandavirus} in the family \textit{Phenuiviridae}. Infected patients experience high fever, gastrointestinal symptoms (vomiting and diarrhea), severe reductions in platelet count and white blood cell count, and elevated liver enzymes indicating liver damage. In severe cases, the disease can lead to failure of multiple organs, with a fatality rate of $12-50\%$ \cite{he2020severe}. This high mortality rate makes SFTS a significant public health concern in endemic regions. SFTSV can infect various animal species including rodents, dogs, and goats, which typically show no symptoms. These asymptomatic animals serve as viral reservoirs that maintain circulation in natural ecosystems. Humans become infected primarily through tick bites and human-to-human transmission can occur through direct contact with infected blood or body fluids~\cite{wu2020transmission}.

\subsection{Data collection}
This study uses environmental and epidemiological data related to SFTS in Dalian, Liaoning Province, China. Monthly confirmed human SFTS cases from 2011 to 2022 were obtained from the Liaoning Province Center for Disease Control and Prevention, providing a community-level view of transmission dynamics. Monthly mean temperature over the same period was extracted from the National Climatic Data Center (NCDC) Integrated Surface Database (\href{https://ftp.ncdc.noaa.gov/pub/data/isd-lite/}{data}).

\begin{figure}[h!]
    \centering
    \includegraphics[width=0.8\linewidth]{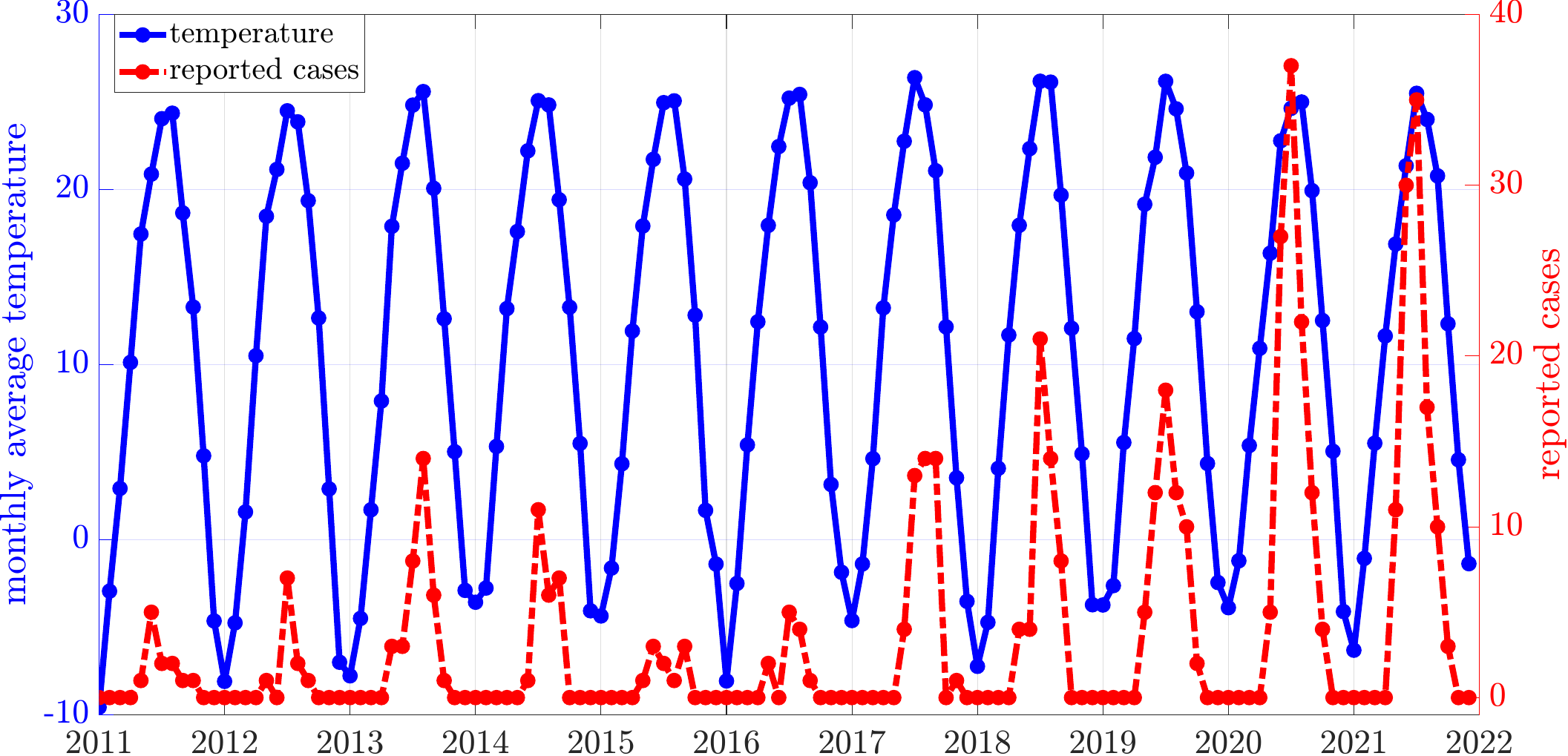}
    \caption{Temporal connection between environmental and epidemiological factors in Dalian from 2011 to 2022 is illustrated. The blue solid line depicts the monthly average ambient temperature, whereas the red dashed line shows the monthly count of confirmed SFTS cases reported.}
    \label{fig:sfts_data}
\end{figure}

Figure~\ref{fig:sfts_data} compares seasonal patterns of temperature and SFTS incidence, showing a clear visual association that supports the hypothesis that climatic factors, particularly temperature, influence tick activity and human exposure. The local ecological setting in Dalian, characterized by a high abundance of \textit{Haemaphysalis longicornis} ticks and typical host species (rodents for nymphs and goats for adults), makes this region well suited for studying tick-borne transmission dynamics \cite{CSIAM-LS-1-2}.

All datasets were aligned to ensure that epidemiological and environmental variables were temporally synchronized for model fitting~\cite{CSIAM-LS-1-2}. This curated dataset enables the construction of high-dimensional, temporally aligned feature matrices suitable for training and validating SINDy-based models.

\subsection{SINDy-based distributed delay model discovery}\label{sec:sfts_sindy}

\subsubsection{Mechanistic baseline model}\label{sec:ode_tick_model}

The reference paper~\cite{CSIAM-LS-1-2} develops a mechanistic, stage-structured compartmental model formulated as a system of ODEs to understand SFTS transmission dynamics in Dalian. The ODE system comprises sixteen state variables representing tick populations at different life stages and infection classes, plus rodent and goat host populations, yielding a coupled system of twenty differential equations for the complete tick-host-virus dynamics. This model explicitly incorporates the complete life cycle of the primary vector, {\it Haemaphysalis longicornis}, tick populations across eight distinct life stages: eggs, questing larvae, feeding larvae, questing nymphs, feeding nymphs, questing adults, feeding adults, and egg-laying adults. Each of these life stages is further subdivided into subcategories (susceptible or infected), recognizing that different developmental stages exhibit distinct infection dynamics, host-seeking behavior, and transmissibility. This age-structured approach is biologically motivated by the fact that immature ticks (larvae and nymphs) predominantly feed on rodents, while adult ticks preferentially parasitize larger hosts such as goats, creating distinct ecological roles and transmission opportunities that a simplified, age-structured model would necessarily overlook.

A key feature of the mechanistic model is its explicit treatment of three coexisting transmission pathways through which SFTSV circulates within the tick-host system. The model incorporates systemic transmission that occurs when susceptible ticks acquire infection while feeding on infected animal hosts or when infected ticks infect susceptible hosts during feeding. Co-feeding transmission, also included in the model, is a non-systemic route wherein susceptible and infected ticks simultaneously feed on the same host, enabling virus transmission through host-mediated mechanisms. The model formulation also incorporates transovarial transmission that represents vertical transmission from infected adult female ticks to their offspring through eggs, providing a mechanism for virus persistence even when host infections are low.

Temperature dependence is encoded by making most biological rates functions of ambient temperature: tick development rates, mortality rates at each life stage, and host-seeking activity all follow temperature-dependent functional forms fitted to experimental data~\cite{CSIAM-LS-1-2}. This is crucial because the tick life cycle shapes both the size and timing of infected tick populations available to transmit the virus to humans. The model also includes logistic birth dynamics for rodent and goat populations, natural mortality at each stage, and a one month reporting delay to represent the lag between tick exposure and confirmed human SFTS diagnosis.

The crucial quantity connecting tick host dynamics to observable human cases is the monthly case incidence, which is computed through an integral equation that aggregates the contributions of infected questing nymphs and infected questing adults over a one month surveillance period, weighted by their respective temperature modulated attachment rates to humans. Specifically, the predicted number of newly reported human SFTS cases during month $t$ is expressed as

\begin{equation}\label{eq:reference_model_cases}
\hat{y}(t)=\int_{0}^{\sigma}(\Lambda_{N}(t-a){p}_{N}(t-a)N_{qi}(t-a)+\Lambda_{A}(t-a){p}_{A}(t-a)A_{qi}(t-a)) \dd a,
\end{equation}
where $\sigma$ denotes the surveillance interval, $N_{qi}$ and $A_{qi}$ represent the densities of infected questing nymphs and infected questing adults at time $t-a$, respectively, and $\Lambda_N$ and $\Lambda_A$ are temperature dependent host attachment rates for nymphs and adults. These attachment rates are formulated as $\Lambda_{\beta} = \alpha_{\beta} e^{\omega T}$, where $\alpha_{\beta}$ represents temperature-dependent strength coefficients, $\omega$ encodes the strength of human outdoor activity and $T$ denotes ambient temperature. This formulation represents the key biological insight that human case incidence depends not merely on the current state of the system but rather on the accumulated history of infected tick populations weighted by their behavioral activity, a concept intimately connected to the distributed delay framework employed in the subsequent SINDy based analysis.

\subsubsection{SINDy-based model discovery}\label{sec:sindy_sfts_dalian}

This section demonstrates how SINDy framework in Section~\ref{sec:dde_re} is used to identify the kernel of a distributed delay model for SFTS transmission in Dalian. Our aim is to infer the underlying transmission dynamics and to forecast future SFTS incidence from multivariate observational data. We employ a data-driven approach using SINDy to identify sparse dynamical relationships directly from observational time series, complementing the mechanistic modelling of transmission dynamics in the tick-host ecosystem. The dataset includes four main variables: monthly human case counts, environmental temperature, and abundance estimates for two tick life stages (nymphs and adults). Tick abundance time series are derived from the mechanistic transmission model in~\cite{CSIAM-LS-1-2}, providing synthetic yet ecologically realistic population estimates. By combining these data sources within SINDy framework, we capture interactions between tick life-stage dynamics, environmental conditions, and human incidence, and recover a sparse representation of the underlying distributed delay kernel suitable for forecasting.

The monthly records of SFTS cases from Dalian spanning 2011 to 2022, along with daily temperature data, serve as input for model identification. The data is divided into training and validation subsets, with the first eight years (96 data points) dedicated to model training and parameter estimation, while the remaining years are set aside for validation and forecasting evaluation. Polynomial SINDy libraries were constructed in combination with trapezoidal quadrature $K=100$ for all the experiments.
In this study, we set $\lambda$ to be very low ($\lambda=10^{-10}$) to ensure the coefficients match with the values found in the reference model \cite{CSIAM-LS-1-2}. This choice prioritizes biological interpretability over maximum sparsity, avoiding the removal of terms just to obtain a more sparse model.

\subsubsection{Model performance: incorporating incidence, outdoor activities and temperature only}

A library of potential kernel functions is created, beginning with terms that include human cases (denoted as $C$), temperature data (denoted as $T$) and their combinations~\eqref{eq:sindy_dd}. Next, an exponential term representing the impact of outdoor activity, expressed as $e^{\omega T}$, is added to the database~\eqref{eq:sindy_dd_exp}. Additionally, a model to account for the distributed memory effect is constructed using numerical quadrature, where the time-delay interval is discretized to establish the regression matrix, and sparse regression is employed to identify the minimal set of functional terms needed to reflect the time-series behavior observed in prior cases.

We investigate two related model forms. The baseline model is
\begin{equation}\label{eq:sindy_dd}
\hat{y}(t)=\int_{0}^{\sigma}f(C(t-a), T(t-a)) \dd a,
\end{equation}
whilst the second incorporates the outdoor activity term
\begin{equation}\label{eq:sindy_dd_exp}
\hat{y}(t)=\int_{0}^{\sigma}f(C(t-a), T(t-a),e^{\omega T(t-a)}) \dd a,
\end{equation}
where $\hat{y}$ denotes the newly reported cases of human SFTS during a specific time period $\sigma$, with the reporting period set to one month. The parameter $\omega$ is externally optimized by applying particle swarm optimization (PS). Data splitting is conducted to assess the generalizability of the model. 

\begin{figure}[h!]
\centering
\includegraphics[width=0.8\linewidth]{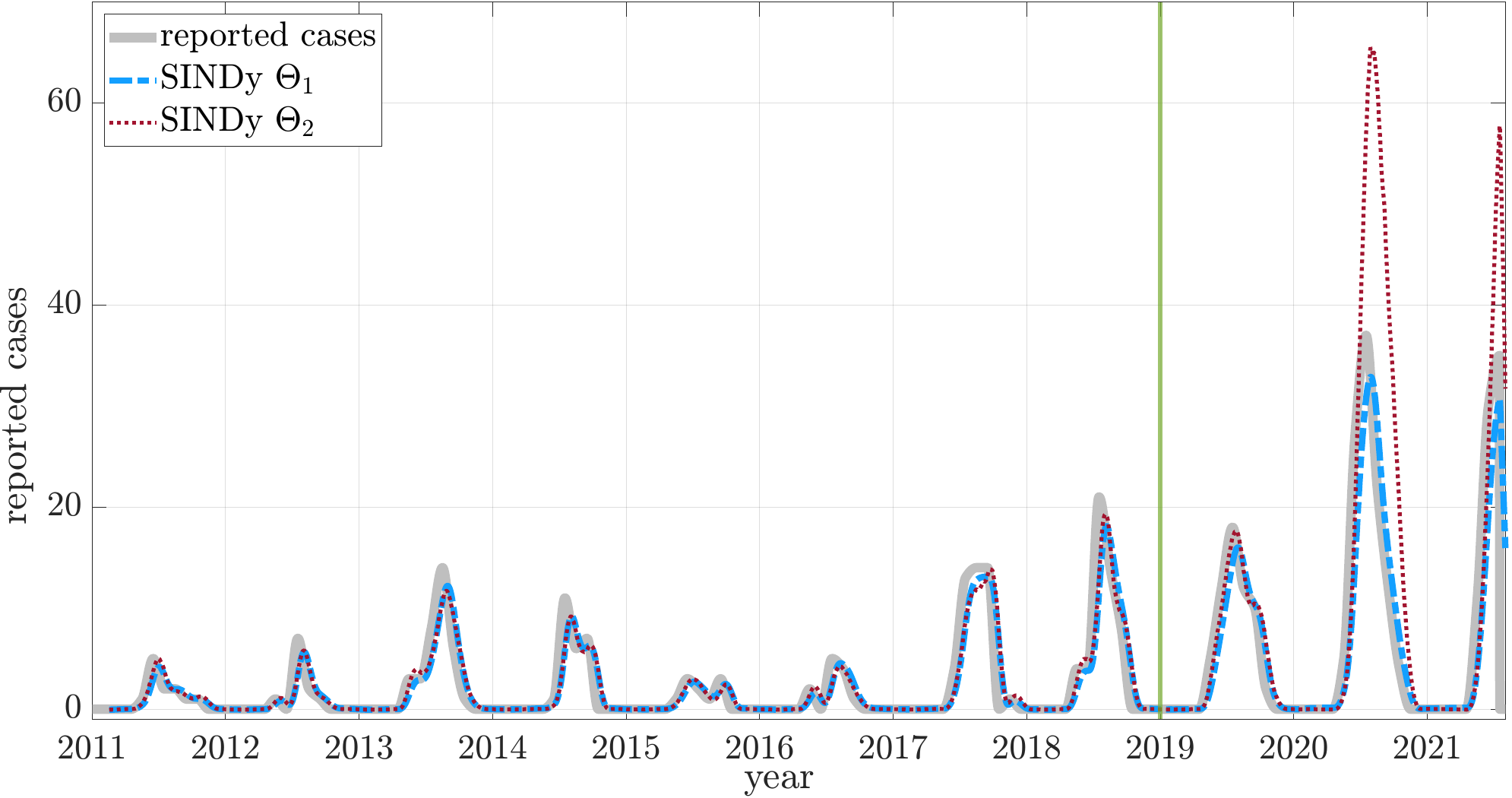}
\caption{Predictive performance of SINDy distributed delay model in forecasting SFTS cases involved utilizing two distinct polynomial order libraries within the candidate function set \eqref{eq:sindy_dd}. The comparison of actual human case numbers (solid grey line) against SINDy predictions for polynomial degree one (blue line) and polynomial degree two (red line).}\label{fig:sindy_sfts1}
\end{figure}

Figure \ref{fig:sindy_sfts1} shows a comparison between the actual reported human cases and the predictions made by SINDy model, which employs two different libraries with varying polynomial orders using human case data and temperature data \eqref{eq:sindy_dd}. The figure demonstrates the capacity of the model to capture the temporal dynamics associated with the incidence of SFTS, with predictions reasonably approximating the observed case counts during both the training and validation phases. The visual comparison indicates that both polynomial orders provide reasonable approximations, although performance differences are noticeable across various time intervals.

\begin{figure}[h!]
\centering
\includegraphics[width=0.8\linewidth]{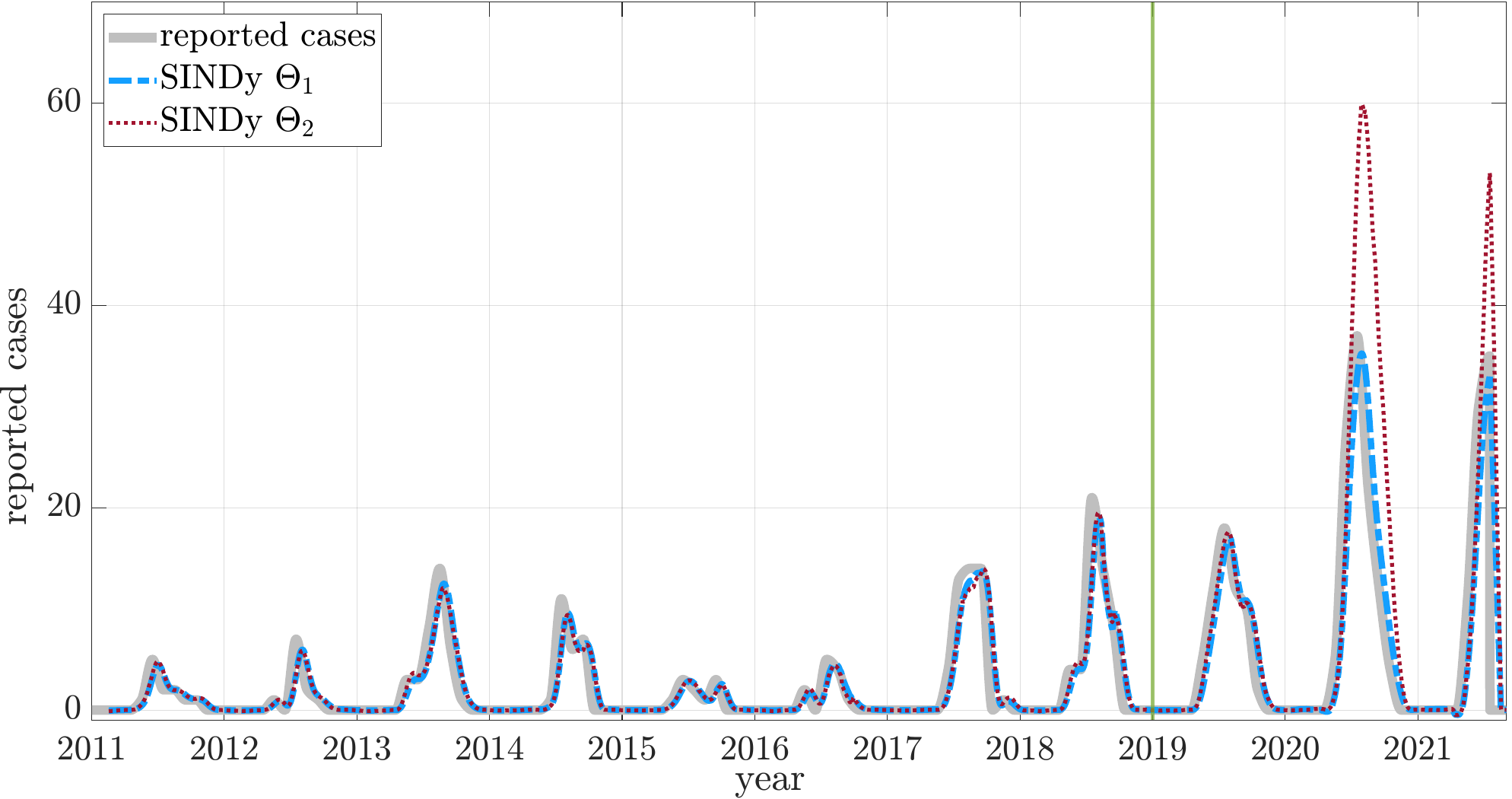}
\caption{Comparison of actual human cases of SFTS (grey line) with SINDy distributed delay model predictions using polynomial degrees one (blue) and two (red). The candidate library includes an exponential term $e^{\omega T}$ \eqref{eq:sindy_dd_exp} to represent temperature-dependent human outdoor activity effects. Polynomial degree one yields improved predictive accuracy.}\label{fig:sindy_sfts1_exp}
\end{figure}

Figure \ref{fig:sindy_sfts1_exp} illustrates a comparison between observed human case data and the predictions made by SINDy model when an exponential term \eqref{eq:sindy_dd_exp} is incorporated into the library for different polynomial orders. This exponential term accounts for the influence of temperature on (human) outdoor activities. The findings indicate that adding this term results in predictive performance similar to the baseline model \eqref{eq:sindy_dd}, with a polynomial order of one showing significantly better alignment between predictions and observed cases than higher-order polynomial terms.

\begin{table}[htbp]
\centering
\begin{tabular}{l cc cc}
\toprule
& \multicolumn{2}{c}{$\xi$ (degree = 1)} & \multicolumn{2}{c}{$\xi$ (degree = 2)} \\
\cmidrule(lr){2-3} \cmidrule(lr){4-5}
{\textit{library terms}} & \eqref{eq:sindy_dd} & \eqref{eq:sindy_dd_exp} &  \eqref{eq:sindy_dd} &  \eqref{eq:sindy_dd_exp} \\
\midrule
$C$ & $3.04\times 10^{-2}$ & $3.26\times 10^{-2}$ & $6.95\times 10^{-2}$ & $6.20\times 10^{-2}$ \\ 
$T$ & $1.77\times 10^{-4}$ & $1.66\times 10^{-4}$ & $1.35\times 10^{-4}$ & $2.49\times 10^{-4}$ \\
$e^{\omega T}$ & - & - & - & $1.85\times 10^{-9}$ \\
$C^2$ & - & - & $2.33\times 10^{-2}$ & $2.09\times 10^{-4}$ \\
$C \cdot T$ & - & - & $-1.66\times 10^{-3}$ & $-1.35\times 10^{-3}$ \\
$T^2$ & - & - & $-7.48\times 10^{-6}$ & $-1.57\times 10^{-5}$ \\
\midrule
\multicolumn{5}{l}{\textit{RMSE}} \\
training & $1.67\times 10^{0}$ & $1.63\times 10^{0}$ & $1.62\times 10^{0}$ & $1.61\times 10^{0}$ \\
validation & $4.96\times 10^{0}$ & $4.99\times 10^{0}$ & $1.20\times 10^{1}$ & $1.04\times 10^{1}$\\
\bottomrule
\end{tabular}
\caption{Sparse coefficient matrix for the distributed delay models, as defined in \eqref{eq:sindy_dd} and  \eqref{eq:sindy_dd_exp}, is identified and estimated using the STLS regression approach with polynomial degrees one and two. The table presents the nonzero coefficients related to the state variables $C$ (human cases), $T$ (temperature), $e^{\omega T}$ (exponential term representing human outdoor activity), along with their nonlinear interactions. RMSE values are provided for both the training and validation datasets.}
\label{tab:comparison1}
\end{table}

\subsubsection{Model performance: incorporating incidence, outdoor activities and temperature into disease transmission dynamics}

Expanding on the results of the previous analysis, we extend the model by integrating data of infected tick populations with human case records and temperature data. This improved method takes advantage of the multivariate nature of the transmission system, acknowledging that the dynamics of SFTS are closely linked to the presence of infected vectors at various life stages. Data on infected nymph and adult tick populations are obtained from simulations of the mechanistic transmission model detailed in \cite{CSIAM-LS-1-2}. As in the previous section, we use the same dataset temporal division. By incorporating these complementary data streams into SINDy framework, the approach captures the complex interdependencies between tick population dynamics, environmental factors, and human infection rates, thus allowing the identification of a more comprehensive representation of the underlying distributed delay kernel. This relationship can be expressed through the following formulation

\begin{equation}\label{eq:sindy_dd_exp_tina}
\hat{y}(t)=\int_{0}^{\sigma}f( C(t-a), T(t-a),e^{\omega T(t-a)},I_{N_{q}}(t-a),I_{A_{q}}(t-a)) \dd a
\end{equation}
where $I_{N_{q}}$ and $I_{A_{q}}$ represent infected questing nymph and adult tick populations, respectively.

\begin{figure}[h!]
\centering
\includegraphics[width=0.8\linewidth]{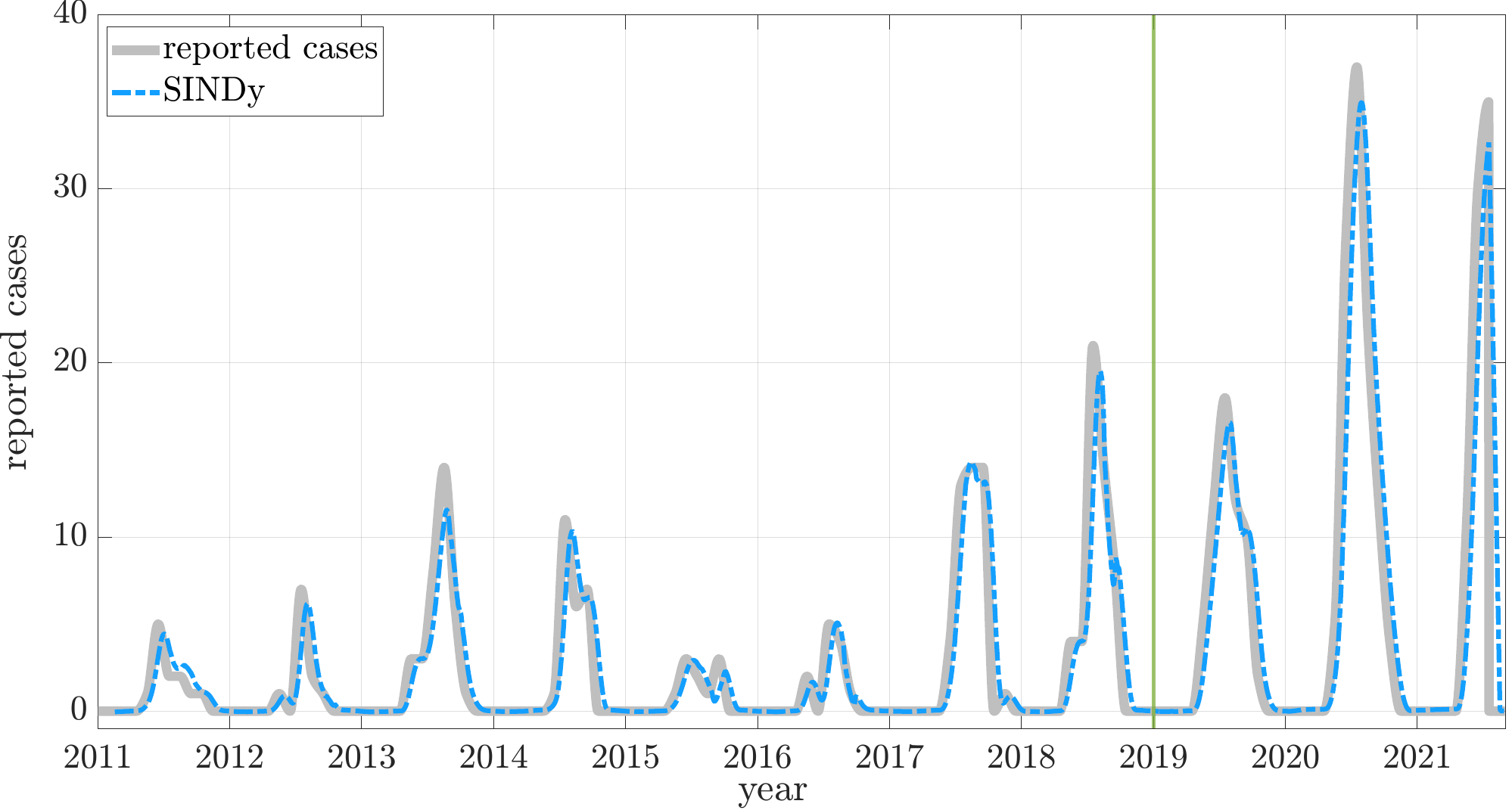}
\caption{Validation of the enhanced SINDy model that integrates data of infected tick populations. The candidate library includes infected nymph and adult populations, as well as environmental and human epidemiological factors \eqref{eq:sindy_dd_exp_tina}. Solid grey line: recorded SFTS cases; blue line: SINDy forecasts.}\label{fig:sfts_iae}
\end{figure}

Figure \ref{fig:sfts_iae} compares actual human cases with SINDy predictions using data that include human cases, temperature, and infected nymph and adult tick populations in the library. The inclusion of these data appears to capture additional dynamics not represented by environmental and epidemiological variables alone.

\begin{table}[htbp]
\centering
\small 
\begin{tabular}{lccccc}
\toprule 
{\textit{library terms}}& $C$ & $T$ & $e^{\omega T}$ & $I_{N_{q}}$ & $I_{A_{q}}$ \\
$\xi$ & $2.61\times 10^{-2}$ & $-7.75\times 10^{-4}$ & $3.55\times 10^{-8}$ & $6.58\times 10^{-5}$ & $-7.03\times 10^{-8}$ \\
\toprule
{\textit{RMSE}} & \multicolumn{2}{c}{training} & \multicolumn{3}{c}{validation} \\
& \multicolumn{2}{c}{$1.51\times 10^{0}$} & \multicolumn{3}{c}{$4.88\times 10^{0}$} \\
\bottomrule

\end{tabular}
\caption{Identified sparse coefficient matrix for the distributed delay model \eqref{eq:sindy_dd_exp_tina} that includes human cases ($C$), temperature ($T$), exponential outdoor activity ($e^{\omega T}$), populations of infected questing nymphs ($I_{N_{q}}$), and populations of infected questing adult ticks ($I_{A_{q}}$) as state variables. Coefficients are recovered using STLS regression at a polynomial degree of one. The table presents RMSE values for both training and validation phases.}
\label{tab:comparison2}
\end{table}


The final identified model for predicting cases of SFTS has an integral equation form, where new cases are expressed through a sparse linear combination of terms of the kernel. The structure of this model has a limited number of non-zero coefficients, fitting well into the objective of sparsity promotion of SINDy approach. The distributed delay model developed in SINDy approach has the best performance in fitting the training dataset, as evidenced by the reasonably small RMSE values reported. The RMSEs reported here are quite reasonable for the noisy surveillance datasets when taking into consideration the fact that there are only 96 monthly samples used for training. 
SINDy models are able to reproduce the seasonality of the epidemics (\Cref{fig:sindy_sfts1,fig:sindy_sfts1_exp,fig:sfts_iae}).

\begin{figure}[h!]
\centering
\includegraphics[width=0.8\linewidth]{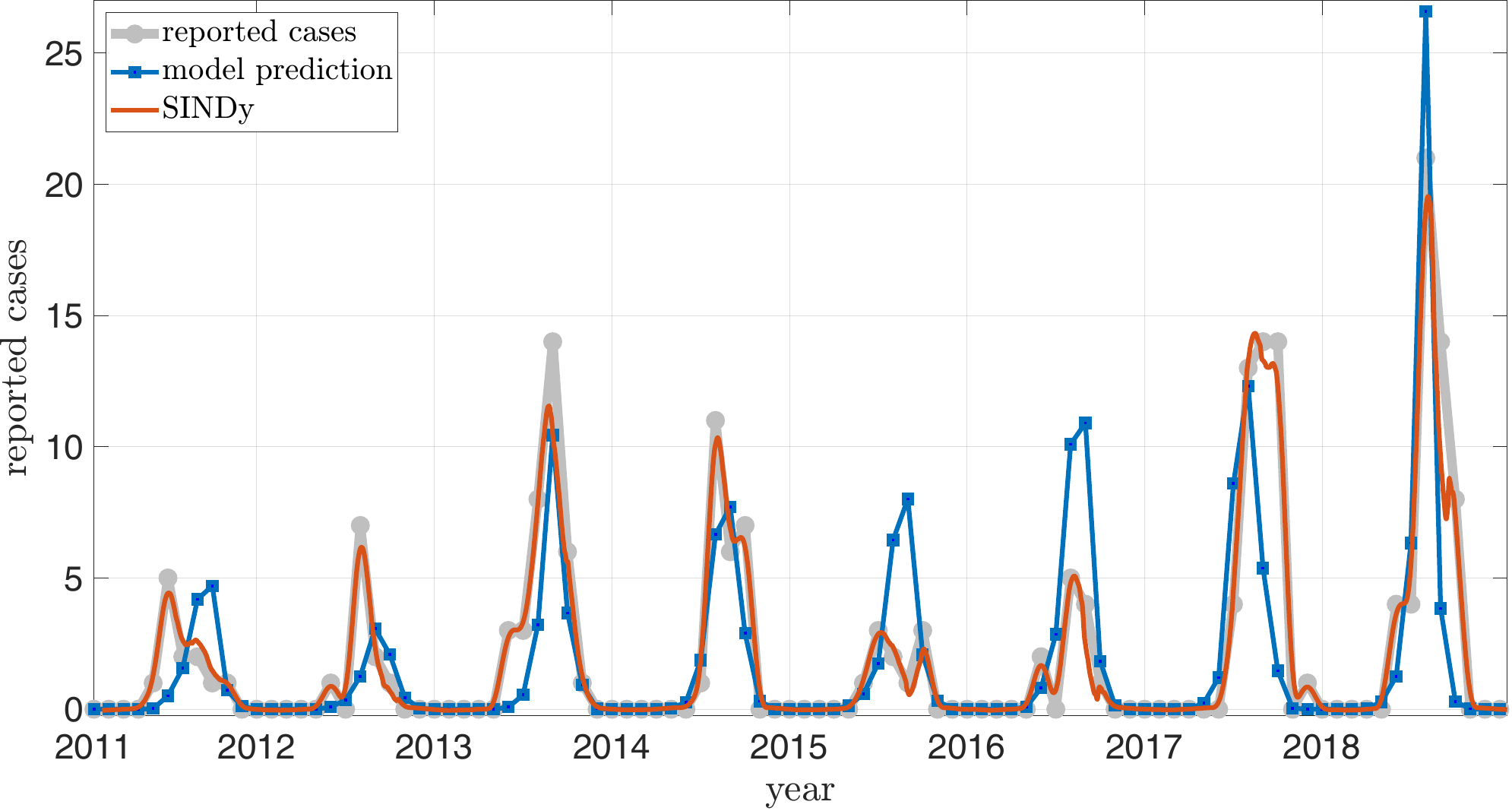}
\caption{Comparison of monthly reported cases (grey) with model predictions (blue) \cite{CSIAM-LS-1-2} and SINDy-based predictions (red) \eqref{eq:sindy_dd_exp_tina}. }\label{fig:sfts_ode_sindy}
\end{figure}

Figure \ref{fig:sfts_ode_sindy} presents the comparison of the trajectories produced by the mechanistic model \cite{CSIAM-LS-1-2} and SINDy discovered model \eqref{eq:sindy_dd_exp_tina} with real monthly reported human cases data. 

\subsection{Parameter sensitivity analysis}
Sensitivity analysis plays a key role in both validating and interpreting data-driven epidemiological models. It shows how model predictions change when parameters are perturbed, which helps assess whether the identified parameters represent real characteristics of the transmission dynamics or numerical artifacts of the optimization procedure. It also improves interpretability by revealing which parameters most strongly influence predictions, thereby highlighting the primary processes in SFTS transmission in Dalian.

SINDy distributed delay model involves two parameters that require sensitivity analysis. The outdoor activity parameter $\omega$ (in the exponential term ${e}^{\omega T}$) captures how temperature affects human exposure to infected ticks through changes in outdoor behavior. The integration window $\sigma$ (defining the reporting period) specifies the time span over which past environmental conditions and vector populations contribute to current human infections.

From an epidemiological  perspective, the model performance and model-based prediction are expected to be sensitive to  both parameters. The parameter $\omega$ should lie within a realistic range: extreme values would imply either negligible or implausibly strong temperature effects on human exposure. The parameter $\sigma$ should be consistent with the characteristic timescales of the tick life cycle and infection process in the study area; substantial deviations from this window would be expected to reduce predictive accuracy if the identified value truly reflects the intrinsic memory of the system rather than a numerical error.

In this section, we systematically adjust these parameters within biologically informed ranges and measure the resulting changes in validation error metrics.

\subsubsection{Outdoor activity parameter $\omega$}

The exponential human activity $e^{\omega T}$ modulates the temperature dependent contact rate between humans and infected tick populations. This term captures the relationship between temperature and outdoor activities that increase human exposure to tick-borne pathogens.

We evaluate the performance of the model over a range of $\omega$ values.
\begin{equation}
\omega \in \left[0.1 \cdot \omega_{\text{opt}}, 1 \right]
\end{equation}
where $\omega_{\text{opt}}$ denotes the optimal value obtained through PS during the training phase. For each candidate value $\omega_i$ on a discretized grid of 15 points, we reconstruct SINDy model using Algorithm~\ref{alg:distributed-delay-sindy} and compute the validation metrics. The corresponding validation error is computed as

\begin{equation}
\begin{aligned}
\epsilon(\omega_i) &= \sqrt{\frac{1}{n} \sum_{j=1}^{n} \left(y_j - \hat{y}_j(\omega_i)\right)^2}.
\end{aligned}
\end{equation}

\subsubsection{Integration window parameter $\sigma$}

The parameter $\sigma$ defines the temporal extent of the distributed delay kernel, representing the epidemiologically relevant time window over which past environmental conditions and tick populations influence the current incidence of SFTS. We extend the analysis range compared to previous studies to explore potential long-term memory effects:
\begin{equation}
\sigma \in [0.25, 3.0] \text{ months}.
\end{equation}
For each $\sigma_k$, the distributed delay integral is approximated via quadrature method (i.e., trapezoidal).
\begin{equation}
\hat{y}(t; \sigma_k) = \int_0^{\sigma_k} f(C(t-a), T(t-a),e^{\omega T(t-a)},I_{N_{q}}(t-a),I_{A_{q}}(t-a)) \dd a.
\end{equation}


\begin{figure}[h!]
\centering
\includegraphics[width=0.8\linewidth]{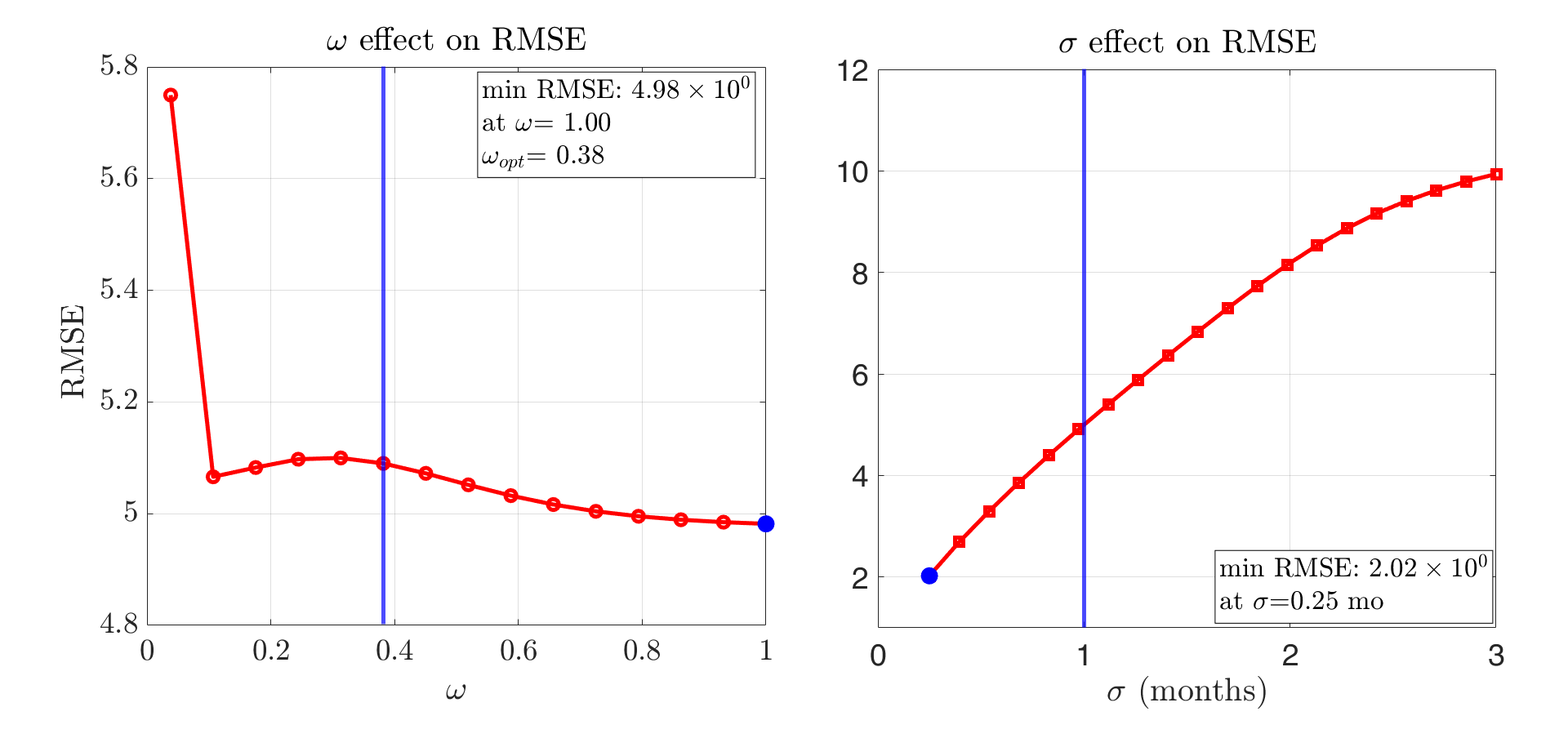}
\caption{Sensitivity analysis was conducted to evaluate how resilient the model predictions are to changes in two critical parameters \eqref{eq:sindy_dd_exp}: the outdoor activity coefficient $\omega$ and the integration window parameter for the human case reporting period $\sigma$, while also including the validation RMSE across a discretized range of parameter values,  vertical blue line in the $\omega$ panel represents the optimal value obtained during training phase and vertical blue line in $\sigma$ panel represents the reference one-month human cases reporting period.}\label{fig:senst_exp_omega_sigma}
\end{figure}

\begin{figure}[h!]
\centering
\includegraphics[width=0.8\linewidth]{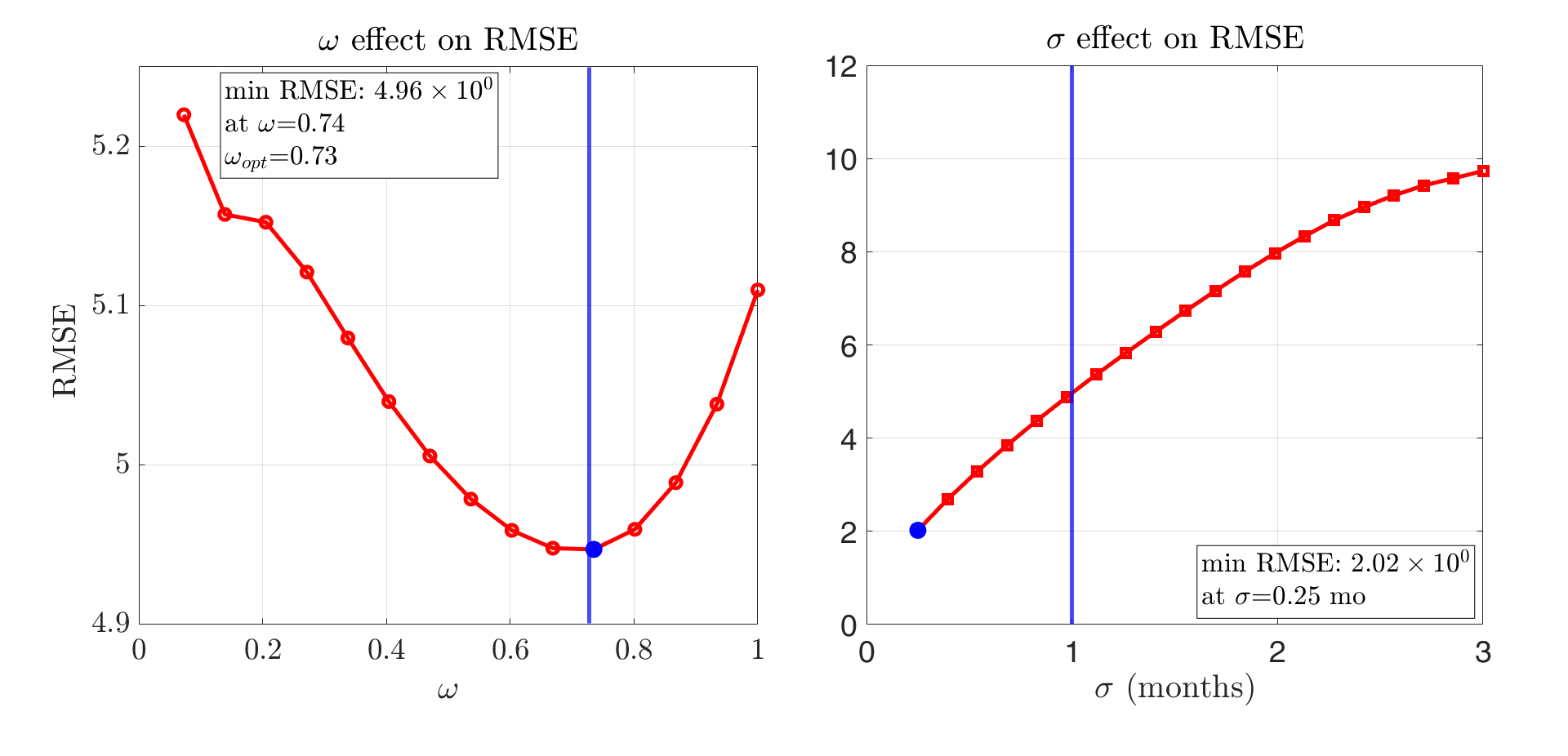}
\caption{Extended parameter sensitivity analysis for the comprehensive distributed delay model incorporates tick population data, which includes both infected nymphs and adults, along with human cases and temperature \eqref{eq:sindy_dd_exp_tina}. This analysis also includes the validation RMSE across different combinations of the outdoor activity parameter $\omega$ and the integration window parameter for the human case reporting period $\sigma$, vertical blue line in the $\omega$ panel represents the optimal value obtained during training phase and vertical blue line in $\sigma$ panel represents the reference one-month human cases reporting period. }\label{fig:senst_sfts_iae}
\end{figure}

Figure \ref{fig:senst_exp_omega_sigma} presents a parameter sensitivity analysis of the outdoor activity parameter $\omega$ and the human case reporting period $\sigma$ for model \eqref{eq:sindy_dd_exp}. The sensitivity analysis provides a quantitative assessment of how variations in these parameters influence the performance of the model and the accuracy of the prediction. 

Figure \ref{fig:senst_sfts_iae} displays the parameter sensitivity analysis of the outdoor activity parameter and the human case reporting period while employing data on human cases, temperature, and infected nymph and adult tick populations for model \eqref{eq:sindy_dd_exp_tina}. The sensitivity profiles provide a comprehensive assessment of the influence of the parameters across the expanded model space, allowing the identification of parameters that have the greatest impact on prediction accuracy.

\section{Conclusions}\label{sec:conclusion}

Predicting the human burden of vector-borne diseases in a changing environment from limited time series data—such as human incidence and temperature records—requires identification of hidden mechanisms underlying disease transmission in the vector–host zoonotic cycle and vector–human contact processes that generate the observed human incidence patterns. These mechanisms are typically highly nonlinear, arising from complex interactions among vectors, hosts, and humans. Moreover, they often incorporate memory effects due to delays in pathogen transmission, vector development, and human case reporting.

The Sparse Identification of Nonlinear Dynamics (SINDy) framework, pioneered in \cite{bpk16}, provides a powerful data-driven approach for uncovering such hidden mechanisms from time series data and has become an important tool in the rapidly expanding field of scientific discovery driven by data.
Unlike previous extensions of the SINDy framework,  our study explicitly targets epidemiological systems with delayed transmission and couples data-driven discovery with mechanistic subsystems.
In particular, we demonstrate that an appropriate extension of the SINDy framework to systems with distributed memory effects [6] can effectively uncover transmission dynamics of severe fever with thrombocytopenia syndrome (SFTS), a tick-borne disease, using limited human incidence and local temperature time series from the study region (Dalian, Liaoning Province, China), without imposing predefined mechanistic assumptions on the human case reporting process.

Our simulations further show that this data-driven sparse identification approach is most effective for predicting human incidence when it is combined with mechanistic modeling of infection transmission in the tick–host subsystem. Incorporating known transmission pathways \cite{CSIAM-LS-1-2}, informed by laboratory experiments, surveillance data, and previous investigations, significantly enhances predictive performance for the highly nonlinear tick–host–human interaction system. This hybrid strategy leverages the complementary strengths of data-driven discovery and mechanistic understanding.

The SINDy framework adapted for distributed memory effects assumes that memory spans a finite time interval and that the kernel function can be expressed as a sparse linear combination of basis functions from a candidate library. Both the integration interval and selected candidate functions may involve parameters that must be optimized through error-based estimation. In our case study, these include the exponent governing human activity effects and the temporal extent of the distributed memory kernel. We show that the proposed SINDy-based framework facilitates systematic sensitivity analysis, yielding comprehensive profiles that quantify the influence of key parameters across the expanded model space and identify those with the greatest impact on predictive accuracy.
We should emphasize that our framework is readily transferable to other vector-borne systems characterized by nonlinear interactions and delayed transmission.

While this data-driven approach prioritizes predictive performance, it may trade mechanistic transparency for computational efficiency and interpretability. As a result, inferred models are sometimes best compared through their dynamical trajectories rather than through explicit mechanistic interpretation. 
While not fully mechanistic, the identified integral terms correspond to biologically and behaviorally meaningful delayed processes, and the identified sparse integral representations offer a practical and flexible framework for epidemiological forecasting and public health resource allocation during peak transmission periods. Although we address the issue of uncertainty in our sensitivity analysis, future work should focus on incorporating uncertainty quantification for probabilistic forecasting, integrating additional environmental drivers, and developing multi-regional models applicable to endemic settings.

\appendix
\section{Algorithms}\label{sec:appendix}


Algorithm \ref{alg:distributed-delay-sindy} summarizes SINDy procedure for reconstructing REs of the form given in equation \eqref{eq:re} from data, as presented in Section \ref{sec:dde_re} using the quadrature-based method. In this approach, the delay is discretized and based on the epidemiological and environmental variables shifted in time along with associated weights, the sparse regression is performed on these to generate the coefficient matrix $\Xi$, indicating the distributed delay function.

In this model, the exponential term $e^{\omega T}$ depends on the human outdoor activity $\omega$, of which the value is not known a priori. To approximate the values of $\omega$, Algorithm \ref{alg:distributed-delay-sindy} is placed inside another routine, Algorithm \ref{alg:adaptive-distributed-sindy}. These candidate values of $\omega$ are then used to approximate SINDy reconstruction error. We define this error for REs as

\begin{equation}\label{eq:optim_error}
\varepsilon({\omega})\coloneqq\|\mathbf{X}-\Theta(\omega)\Xi(\omega)\|_{2},
\end{equation}
where $\Theta(\omega)$ and $\Xi(\omega)$ obtained for that specific $\omega$. An external optimizer like PS \citep{bonyadi2017,kennedy1995,shi1998} is employed over $\omega$ to optimize $\varepsilon(\omega)$ and obtain an optimal parameter $\omega^{*}$ along with a corresponding sparse kernel $\Xi^{*}$. This two-level strategy helps to obtain both a distributed delay kernel as well as an expression for human outdoor activity in a unified manner.

\begin{algorithm}[H]
\caption{SINDy for Distributed Delay and Renewal equations}
\label{alg:distributed-delay-sindy}
\begin{algorithmic}[1]
\State \textbf{Input:} Data: $\mathbf{C}$, $\mathbf{T}$, $\mathbf{I}_{N_q}$, $\mathbf{I}_{A_q}$, function library $ \mathcal{F}$, delay interval $[0, s]$, quadrature rule parameters, regularization parameter $\lambda$
\State \textbf{Output:} Sparse coefficient matrix $\Xi$ representing kernel $$g(a, C(t-a), T(t-a), I_{N_q}(t-a), I_{A_q}(t-a))$$

\State \textbf{Step 1: Quadrature Setup}
\State Discretize the integration interval $[0, s]$ into $K$ nodes $s_{k}$ with weights $w_{k}$
\State \quad (e.g., via trapezoidal, Clenshaw-Curtis, or Gaussian quadrature)

\State \textbf{Step 2: Data Preparation}
\For{$k = 1$ to $K$}
    \State Compute delay $s_k$ from integration nodes
    \State Obtain time-shifted data by interpolation:
    \State \quad $\mathbf{C}_{s_{k}} = \mathbf{C}(t - s_k)$, $\mathbf{T}_{s_{k}} = \mathbf{T}(t - s_k)$, $\mathbf{I}_{N_q,s_{k}} = \mathbf{I}_{N_q}(t - s_k)$, $\mathbf{I}_{A_q,s_{k}} = \mathbf{I}_{A_q}(t - s_k)$, $\mathbf{E}_{s_{k}} =e^{\omega \mathbf{T}_{s_{k}}}$
\EndFor

\State \textbf{Step 3: Library Construction}
\State $\Theta_k(a_k) = \big[1, \, s_k, \, \mathbf{C}_{s_{k}}, \, \mathbf{T}_{s_{k}}, \, \mathbf{E}_{s_{k}}, \, \mathbf{I}_{N_q,s_{k}}, \, \mathbf{I}_{A_q,s_{k}}, \, \ldots \big]$
   $\Theta = \sum_{k=1}^{K} w_k \cdot \Theta_k(s_k)$

\State \textbf{Step 4: Sparse Regression}
\For{$j = 1$ to $n$ (each state component)}
    \State Solve sparse regression using STLS or similar:

      \State  $\xi_{j} = \arg\min_{\xi} \left(\left\| \mathbf{X}_{j} - \left( \sum_{k=1}^{K} w_{k} \mathbf{\Theta}_k \right) \xi \right\|_{2}+\lambda\|\xi\|_{1} \right) $
\EndFor
\State Assemble coefficient matrix: $\Xi := (\xi_{1}, \xi_{2},\ldots,\xi_n)$

\State \textbf{Step 5: Model Validation}

\State \textbf{return} $\Xi$, validation metrics
\end{algorithmic}
\end{algorithm}

\begin{algorithm}[H]
\caption{SINDy with Parameter Optimization for Distributed Delays}
\label{alg:adaptive-distributed-sindy}
\begin{algorithmic}[1]
\State \textbf{Input:} Data $\mathbf{C}$, $\mathbf{T}$, $\mathbf{I}_{N_q}$, $\mathbf{I}_{A_q}$, function library $ \mathcal{F}$, parameter search space $\omega$

\State \textbf{Output:} Optimal parameter $\omega^*$ and sparse coefficient vector $\Xi^*$

\Function{$\mathcal{J}$}{$\omega$}
    \State Generate $K$ quadrature nodes $\{s_{k}\}_{k=1}^K$ with weights $\{w_{k}\}_{k=1}^K$ in $[0, s]$ and obtain time-shifted data $\{\mathbf{C}_{a_{k}}, \mathbf{T}_{a_{k}}, \mathbf{I}^q_{N,a_{k}}, \mathbf{I}^q_{A,a_{k}}\}_{k=1}^K$
    
    \State Compute outdoor activity terms: $\mathbf{E}_{s_{k}} = e^{\omega \mathbf{T}_{s_{k}}}$ for all $k$
    \State Construct weighted library matrix $\Theta$ as in Algorithm \ref{alg:distributed-delay-sindy}
    
    \For{$j = 1$ to $n$} Solve sparse regression to get $\xi_{j}(\omega)$ 
    \EndFor
    
    \State $\varepsilon(\omega) = \|X - \Theta(\omega)\Xi(\omega)\|_{2}$

    \State \textbf{return} $\varepsilon(\omega)$
\EndFunction

\State \textbf{External Parameter Optimization} (PS, BO, etc.): $\omega^* = \arg\min_{\omega} \mathcal{J}(\omega)$ \textbf{return} $\omega^*$, $\Xi^*$

\end{algorithmic}
\end{algorithm}


\section*{Acknowledgements}
DB and MT are members of INdAM research group GNCS; DB is a
member of UMI research group ``Modellistica socio-epidemiologica''. JW is a member of the ``Centre of Excellence in Artificial Intelligence for Public Health Advancement'' and a member of the ``Laboratory for Industrial and Applied Mathematics (LIAM)''. The work of DB was partially supported by the Italian Ministry of University and Research (MUR) through the PRIN 2022 project (No. 20229P2HEA) ``Stochastic numerical modelling for sustainable innovation'', Unit of Udine (CUP G53C24000710006). The work of MT was supported by the Italian Ministry of University and Research (MUR) through a PhD grant PNRR DM351/22 (CUP: G23C22001320003). The work of JW was supported in part by the York Research Chair program (grant number: 492108) and by the NSERC-Sanofi Alliance program in Vaccine Mathematics, Modelling, and Manufacturing (517504). The work of XZ was supported by the Natural Science Foundation of China (No. 12571524).
\bibliographystyle{abbrv}
\bibliography{ref}

\end{document}

\endinput